\documentclass[a4paper]{amsart}
\usepackage{graphicx}
\usepackage{amsmath}
\usepackage{amssymb}
\usepackage{stmaryrd}
\usepackage{oldgerm}
\usepackage[all]{xy}

\newcommand{\Hom}{\operatorname{Hom}\nolimits}
\renewcommand{\Im}{\operatorname{Im}\nolimits}

\newcommand{\Pd}{\operatorname{pd}\nolimits}

\newcommand{\depth}{\operatorname{depth}\nolimits}
\newcommand{\Ann}{\operatorname{Ann}\nolimits}

\newcommand{\Tor}{\operatorname{Tor}\nolimits}
\newcommand{\Ext}{\operatorname{Ext}\nolimits}

\newcommand{\m}{\operatorname{\mathfrak{m}}\nolimits}
\newcommand{\az}{\operatorname{\mathfrak{a}}\nolimits}

\newcommand{\cx}{\operatorname{cx}\nolimits}

\newcommand{\V}{\operatorname{V}\nolimits}

\newcommand{\cidim}{\operatorname{CI-dim}\nolimits}
\newcommand{\lcidim}{\operatorname{CI_{_*}-dim}\nolimits}
\newcommand{\gdim}{\operatorname{G-dim}\nolimits}

\newtheorem{theorem}{Theorem}[section]

\newtheorem{corollary}[theorem]{Corollary}

\newtheorem{lemma}[theorem]{Lemma}
\newtheorem{proposition}[theorem]{Proposition}
\newtheorem{definition}[theorem]{Definition}
\theoremstyle{definition}
\newtheorem*{example}{Example}

\theoremstyle{definition}

\theoremstyle{remark}
\newtheorem*{remark}{Remark}
\theoremstyle{definition}

\begin{document}
\title{Complexity test modules}
\author{Petter Andreas Bergh}
\address{Petter Andreas Bergh \newline Mathematical Institute \\ 24-29
  St.\ Giles \\ Oxford \\ OX1 3LB \\ United Kingdom \newline
  \emph{Present address:} \newline Institutt for matematiske fag \\
NTNU \\ N-7491 Trondheim \\ Norway}

\email{bergh@math.ntnu.no}

\thanks{The author acknowledges the support of the Marie Curie
Research Training Network through a LieGrits fellowship} 



\maketitle

\begin{abstract}
A method is provided for computing an upper bound of the complexity of a
module over a local ring, in terms of vanishing of certain cohomology
modules. We then specialize to complete intersections, which are
precisely the rings over which all modules have finite complexity.    
\end{abstract}

\section{Introduction}

The notion of complexity of a module was introduced by Alperin and
Evens in \cite{Alperin}, in order to study modular representations of
finite groups. A decade later, in \cite{Avramov1, Avramov2}, Avramov
introduced this concept for finitely generated modules over local
rings, as a means to distiguish between modules of infinite projective
dimension. The complexity of a module measures the growth of its
minimal free resolution. For example, a module has complexity zero if
and only if its projective dimension is finite, and complexity one
precisely when its minimal free resolution is bounded.

For an arbitrary local ring, not much is known about the modules of
finite complexity. For example, a characterization of these modules in
terms of cohomology does not exist. Even worse, it is unclear whether
every local ring has finite finitistic complexity dimension, that is,
whether there is a bound on the complexities of the modules having
finite complexity. In other words, the status quo for complexity is
entirely different from that of projective dimension; there simply
does not exist any analogue of the Auslander-Buchsbaum formula. The
only class of rings for which these problems are 
solved are the complete intersections; over such a ring \emph{every}
module has finite complexity, and the complexity is bounded by the
codimension of the ring.

In this paper we give a method for computing an upper bound of the
complexity of classes of modules over local rings. This is done by
looking at the vanishing of cohomology with certain ``test'' modules,
using the notion of reducible complexity introduced in
\cite{Bergh1}. We then specialize to the case when the ring is a
complete intersection, and use 
the theory of support varieties, introduced in \cite{Avramov1} and
\cite{AvramovBuchweitz}, both to sharpen our results and to obtain new
ones. 

The paper is organized as follows. In the next section we introduce
the class of modules having ``free reducible complexity'', a slight
generalization of the notion of reducible complexity. We then study
the homological behavior of such modules, proving among other things
some results on the vanishing of (co)homology. In the final section we
prove our main results on complexity testing, specializing at the end
to complete intersections.

\section{Free reducible complexity}

Throughout we let ($A, \m, k$) be a local (meaning also commutative
Noetherian) ring, and we suppose all modules are finitely generated.
For an $A$-module $M$ with minimal free resolution
$$\cdots \to F_2 \to F_1 \to F_0 \to M \to 0,$$ 
the rank of $F_n$, i.e.\ the integer $\dim_k
\Ext_A^n (M,k)$, is the $n$th \emph{Betti number} of $M$, and we
denote this by $\beta_n(M)$. The $i$th \emph{syzygy} of $M$,
denoted $\Omega_A^i (M)$, is the cokernel of the map $F_{i+1} \to
F_i$, and it is unique up to isomorphism. Note that $\Omega_A^0
(M)=M$ and that $\beta_n \left ( \Omega_A^i (M) \right ) =
\beta_{n+i}(M)$ for all $i$. The \emph{complexity} of $M$, denoted
$\cx M$, is defined as
$$\cx M = \inf \{ t \in \mathbb{N} \cup \{ 0
\} \mid \exists a \in \mathbb{R} \text{ such that } \beta_n(M) \leq
an^{t-1} \text{ for } n \gg 0 \}.$$ In general the complexity of a
module may be infinite, in fact the rings for which all modules have
finite complexity are precisely the complete intersections. From the
definition we see that the complexity is zero if and only if the module
has finite projective dimension, and that the modules of complexity
one are those whose minimal free resolutions are bounded,
for example the periodic modules. Moreover, the complexity of $M$
equals that of $\Omega_A^i (M)$ for every $i \geq 0$. 

Let $N$ be an $A$-module, and consider an element $\eta
\in \Ext^t_A(M,N)$. By choosing a map $f_{\eta} \colon
\Omega_A^t(M) \to N$ representing $\eta$, we obtain a
commutative pushout diagram
$$\xymatrix{
0 \ar[r] & \Omega_A^t(M) \ar[r] \ar[d]^{f_{\eta}} &
F_{t-1} \ar[r] \ar[d] & \Omega_A^{t-1}(M) \ar[r]
\ar@{=}[d] & 0 \\
0 \ar[r] & N \ar[r] & K_{\eta} \ar[r] & \Omega_A^{t-1}(M)
\ar[r] & 0 }$$ with exact rows. Note that the module $K_{\eta}$ is
independent, up to isomorphism, of the map $f_{\eta}$ chosen as a
representative for $\eta$. The module $M$
has \emph{reducible complexity} if either the projective dimension
of $M$ is finite, or if the complexity of $M$ is positive and
finite, and there exists a homogeneous element $\eta \in
\Ext_A^*(M,M)$, of positive degree,  such that $\cx K_{\eta} = \cx M-1$,
$\depth K_{\eta} = 
\depth M$, and $K_{\eta}$ has reducible complexity. In this case, the
cohomological element $\eta$ is said to \emph{reduce the complexity}
of $M$. We denote the category of $A$-modules having reducible
complexity by $\mathcal{C}_A^{rc}$.

\sloppy The notion of reducible complexity was introduced in \cite{Bergh1},
where several (co)homology vanishing results for such modules were
given (see also \cite{Bergh3}). These results were previously known to
hold for modules over complete intersections, and so it is not
surprising that over such rings all modules have reducible
complexity. In fact, by \cite[Proposition 2.2]{Bergh1} every module of
finite \emph{complete intersection dimension} has reducible complexity, and
given such a module the reducing process decreases the complexity by
exactly one. Recall that the module $M$ has finite complete
  intersection dimension, written $\cidim_A M < \infty$, if 
there exist local rings $R$ and $Q$ and a diagram $A \to R
\twoheadleftarrow Q$ of local homomorphisms such that $A \to R$ is
faithfully flat, $R \twoheadleftarrow Q$ is surjective with kernel
generated by a regular sequence (such a diagram $A \to R
\twoheadleftarrow Q$ is called a \emph{quasi-deformation} of $A$),
and $\Pd_Q (R \otimes_A M)$ is finite. Such modules were first
studied in \cite{Avramov3}, and the concept generalizes that of
virtual projective dimension defined in \cite{Avramov1}. As the name
suggests, modules having finite complete intersection dimension to a
large extent behave homologically like modules over complete
intersections. Indeed, over a complete intersection ($S,
\mathfrak{n}$) \emph{every} module has finite complete intersection
dimension; the completion $\widehat{S}$ of $S$ with respect to the
$\mathfrak{n}$-adic topology is the residue ring of a regular local
ring $Q$ modulo an ideal generated by a regular sequence, and so $S
\to \widehat{S} \twoheadleftarrow Q$ is a quasi deformation.

It should be commented on the fact that in the original definition of the
notion of reducible complexity in \cite{Bergh1}, the reducing process
did not necessarily reduce the complexity by exactly one. Namely, if
$\eta \in \Ext_A^*(M,M)$ is a cohomological element reducing the
complexity of the module $M$, then the requirement was only that $\cx
K_{\eta}$ be strictly less than $\cx M$, and not necessarily equal to
$\cx M-1$. However, we are unaware of any example where the complexity
actually drops by more than one. On the contrary, there is evidence to
suggest that this cannot happen. If $\cidim_A M$ is finite and 
$$0 \to M \to K \to \Omega_A^n(M) \to 0$$
is an exact sequence, then it is not hard to see that there exists one
``joint'' quasi deformation $A \to R \twoheadleftarrow Q$ such that
all these three modules have finite projective dimension over $Q$ when
forwarded to $R$. In this situation we can apply \cite[Theorem
1.3]{Jorgensen}, a result originally stated for complete
intersections, and obtain a ``principal lifting'' (of $R$) over which
the complexities drop by exactly one. By alterating this procedure, we
can show that the complexity of $K$ cannot drop by more than one. Note
also that this is trivial if $\cx M \le 2$. 

We now define the class of modules we are going to study in this
paper, a class which is a slight generalization of modules of
reducible complexity. Namely, we do not require equality of depth for
the modules involved.

\begin{definition}\label{def}
Denote by $\mathcal{C}_A$ the category of all $A$-modules having finite
complexity. The full subcategory $\mathcal{C}^{frc}_A \subseteq
\mathcal{C}_A$ consisting of the modules having \emph{free reducible
complexity} is defined inductively as follows. 
\begin{enumerate}
\item[(i)] Every module of finite projective dimension belongs to
$\mathcal{C}^{frc}_A$.
\item[(ii)] A module $X \in \mathcal{C}_A$ of positive complexity
belongs to $\mathcal{C}^{frc}_A$ if there exists a homogeneous element
$\eta \in \Ext_A^* (X,X)$ of positive degree such that $\cx
K_{\eta} = \cx X-1$ and $K_{\eta} \in \mathcal{C}^{frc}_A$.
\end{enumerate}
\end{definition}

Thus if the $A$-module $M$ has finite positive complexity $c$, say, then it
has free reducible complexity if and only if the following hold: there
exists a sequence $M=K_0, K_1, \dots, K_c$ of $A$-modules such that for 
each $1 \le i \le c$ there is an integer $n_i \ge 0$ and an
exact sequence
$$0 \to K_{i-1} \to K_i \to \Omega_A^{n_i}(K_{i-1}) \to 0$$
with $\cx K_i =c-i$. Note that the
only difference between this definition and that of ``ordinary'' 
reducible complexity is that we do \emph{not} require that
$\depth K_i$ equals $\depth K_{i-1}$. However, when the ring $A$ is
Cohen-Macaulay, this is always the case (see the remark following
\cite[Definition 2.1]{Bergh1}). Therefore, for such a ring, the
categories $\mathcal{C}^{rc}_A$ and $\mathcal{C}^{frc}_A$ coincide,
that is, a module has reducible complexity if and only if it has free
reducible complexity.

It is quite clear that the inclusion $\mathcal{C}^{rc}_A \subseteq
\mathcal{C}^{frc}_A$ holds, that is, every module of reducible
complexity also has free reducible complexity. In particular, every
$A$-module of finite 
complete intersection dimension belongs to
$\mathcal{C}^{frc}_A$. However, the converse is a priori not true,
that is, a module in $\mathcal{C}^{frc}_A$ need not have reducible
complexity. Moreover, a module having finite complexity need not have
free reducible complexity. We illustrate all this with an example from
\cite{Gasharov}. 

\begin{example}
Let ($A, \m, k$) be the local finite dimensional algebra $k[X_1,
\dots, X_5] / \az$, where $\az \subset k[X_1, \dots, X_5]$ is the
ideal generated by the quadratic forms
\begin{eqnarray*}
& X_1^2, \hspace{3mm} X_2^2, \hspace{3mm} X_5^2, \hspace{3mm} X_3X_4,
\hspace{3mm} X_3X_5, \hspace{3mm} X_4X_5, \hspace{3mm} X_1X_4+X_2X_4 & 
\\
& \alpha X_1X_3+X_2X_3, \hspace{3mm} X_3^2-X_2X_5+ \alpha X_1X_5,
\hspace{3mm} X_4^2-X_2X_5+X_1X_5 &
\end{eqnarray*}
for a nonzero element $\alpha \in k$. By \cite[Proposition
3.1]{Gasharov} this ring is Gorenstein, and the complex
$$\cdots \to A^2 \xrightarrow{d_{n+1}} A^2 \xrightarrow{d_n} A^2 \to
\cdots$$
with maps given by the matrices $d_n = \left ( \begin{smallmatrix} x_1
    & \alpha^n 
x_3+x_4 \\ 0 & x_2 \end{smallmatrix} \right )$ is exact. This sequence
is therefore a minimal free resolution of the module $M:= \Im d_0$,
hence this module has complexity one. If the order of $\alpha$ in $k$
is infinite, then $M$ cannot have 
reducible complexity (recall that the notion of reducible complexity
coincides with that of free reducible complexity since the ring is
Cohen-Macaulay); if $M$ has reducible complexity, then there exists an
exact sequence
$$0 \to M \to K \to \Omega_A^n(M) \to 0$$
in which $K$ has finite projective dimension. As $A$ is
selfinjective, the module $K$ must be free, hence $M$ is a periodic
module, a contradiction. Moreover, if the order of $\alpha$ is finite
but at least $3$, then the argument in \cite[Section 2,
example]{Bergh1} shows that $M$ has reducible complexity but not finite
complete intersection dimension
\end{example}

This example also shows that, in general, a module of finite
\emph{Gorenstein dimension} and finite complexity need not have reducible
complexity. Recall that a module $X$ over a local ring $R$ has finite
Gorenstein dimension, denoted $\gdim_R X < \infty$, if there
exists an exact sequence
$$0 \to G_t \to \cdots \to G_0 \to X \to 0$$
of $R$-modules in which the modules $G_i$ are reflexive and satisfy
$\Ext_R^j(G_i,R) = 0 = \Ext_R^j(\Hom_R(G_i,R),R)$ for $j \ge 1$. Every
module over a Gorenstein ring has finite Gorenstein dimension, in fact
this property characterizes Gorenstein rings. Using this concept,
Gerko introduced in \cite{Gerko}
the notion of \emph{lower complete intersection
  dimension}; the module $X$ has finite lower complete intersection
dimension, written $\lcidim_R X < \infty$, if it has finite Gorenstein
dimension and finite complexity (and in this case $\lcidim_R X =
\gdim_R X$). The Gorenstein dimension, lower complete intersection
dimension, complete intersection dimension and projective dimension of
a module are all related via the inequalities
$$\gdim_R X \le \lcidim_R X \le \cidim_R X \le \Pd_R X.$$
If one of these dimensions happen to be finite, then it is equal to
those to its 
left. Note that the class of modules having (free) reducible complexity and
finite Gorenstein dimension lies properly ``between'' the class of modules
having finite lower complete intersection dimension and the class of modules
having finite complete intersection dimension.  

We now continue investigating the properties of the category of
modules having free reducible complexity. The following result shows that
$\mathcal{C}^{frc}_A$ is closed under 
taking syzygies and preserved under faithfully flat extensions. We omit
the proof because it is analogous to that of \cite[Proposition
2.2]{Bergh1}. 

\begin{proposition}\label{syzygies}
Let $M$ be a module in $\mathcal{C}^{frc}_A$.
\begin{enumerate}
\item[(i)] The kernel of any surjective map $F \twoheadrightarrow M$
  in which $F$ is free also belongs to $\mathcal{C}^{frc}_A$. In
  particular, any syzygy of $M$ belongs to $\mathcal{C}^{frc}_A$. 
\item[(ii)] If $A \to B$ is a faithfully flat local homomorphism, then
  the $B$-module $B \otimes_A M$ belongs to $\mathcal{C}^{frc}_B$.
\end{enumerate}
\end{proposition} 

Note that, in the first part of this result, as opposed to the
corresponding result \cite[Proposition 2.2(ii)]{Bergh1} for modules
belonging to $\mathcal{C}^{rc}_A$, we do not 
need to assume that the ring in question is Cohen-Macaulay. The reason
for this is of course that in the definition of free reducible
complexity, we do not require equality of depth for the modules
involved. By dropping this requirement, the main results from
\cite{Bergh1} on the vanishing of $\Ext$ and $\Tor$ do not carry over
to the category $\mathcal{C}^{frc}_A$ of modules having free reducible
complexity. However, as the following results show, modified versions
of the mentioned results hold for modules belonging to
$\mathcal{C}^{frc}_A$. We prove only the cohomology case; the homology
case is totally similar.  

\begin{proposition}\label{cohomologyvanishing}
Suppose $M$ belongs to $\mathcal{C}^{frc}_A$ and has positive
complexity. Choose modules $M=K_0, K_1, \dots, K_c$, integers $n_1,
\dots, n_c$ and exact sequences 
$$0 \to K_{i-1} \to K_i \to \Omega_A^{n_i}(K_{i-1}) \to 0$$
as in the equivalent definition following Definition
\ref{def}. Then for any $A$-module $N$, the following are equivalent:
\begin{enumerate}
\item[(i)] There exists an integer $t > \max \{ \depth A -
  \depth K_i \}$ such that $\Ext_A^{t+i}(M,N)=0$ for $0 \le i
  \le n_1 + \cdots + n_c$.
\item[(ii)] $\Ext_A^i(M,N)=0$ for $i \gg 0$.
\item[(iii)] $\Ext_A^i(M,N)=0$ for $i > \max \{ \depth A -
  \depth K_i \}$.
\end{enumerate}
\end{proposition} 

\begin{proof}
We have to prove the implication (i) $\Rightarrow$ (iii), and we do
this by induction on the complexity of $M$. If $\cx M=1$, then in the
exact sequence
$$0 \to M \to K_1 \to \Omega_A^{n_1} (M) \to 0$$
the module $K_1$ has finite projective dimension, which by the
Auslander-Buchsbaum formula equals $\depth A - \depth K_1$. Then for
any $i > \max \{ \depth A - \depth K_i \}$ the cohomology group
$\Ext_A^i(K_1,N)$ vanishes, implying an isomorphism 
$$\Ext_A^i(M,N) \simeq \Ext_A^{i+n_1 +1}(M,N).$$
By assumption, the cohomology group $\Ext_A^i(M,N)$ vanishes for $t
\le i \le t + n_1$, and so the isomorphisms just given ensure that
$\Ext_A^i(M,N)=0$ for all $i > \max \{ \depth A - \depth K_i \}$.

Now suppose the complexity of $M$ is at least two. From the assumption
on the vanishing of $\Ext_A^i(M,N)$ and the exact sequence from the
beginning of the proof, we see that $\Ext_A^i(K_1,N)=0$ for
$t \le i \le t+n_2 + \cdots + n_c$. The complexity of 
$K_1$ is one less than that of $M$, hence by induction the
cohomology group $\Ext_A^i(K_1,N)$ vanishes for all
$i > \max \{ \depth A - \depth K_i \}$. The same argument we
used in the case when the complexity of $M$ was one now shows that
$\Ext_A^i(M,N)=0$ for $i > \max \{ \depth A - \depth K_i \}$.
\end{proof}

\begin{proposition}\label{homologyvanishing}
Suppose $M$ belongs to $\mathcal{C}^{frc}_A$ and has positive
complexity. Choose modules $M=K_0, K_1, \dots, K_c$, integers $n_1,
\dots, n_c$ and exact sequences 
$$0 \to K_{i-1} \to K_i \to \Omega_A^{n_i}(K_{i-1}) \to 0$$
as in the equivalent definition following Definition
\ref{def}. Then for any $A$-module $N$, the following are equivalent:
\begin{enumerate}
\item[(i)] There exists an integer $t > \max \{ \depth A -
  \depth K_i \}$ such that $\Tor^A_{t+i}(M,N)=0$ for $0 \le i
  \le n_1 + \cdots + n_c$.
\item[(ii)] $\Tor^A_i(M,N)=0$ for $i \gg 0$.
\item[(iii)] $\Tor^A_i(M,N)=0$ for $i > \max \{ \depth A -
  \depth K_i \}$.
\end{enumerate}
\end{proposition} 

We end this section with a result from \cite{Bergh1}, a result which
will be of 
use in the next section. It shows that when $A$ is Gorenstein, then 
symmetry holds for the vanishing of cohomology between modules of
reducible complexity. This is a generalization of J{\o}rgensen's
result \cite[Theorem 4.1]{Joergensen}, which says that symmetry in the
vanishing of cohomology holds for modules of finite complete
intersection dimension over a local Gorenstein ring.

\begin{proposition}\label{symmetry}\cite[Theorem 3.5]{Bergh1}
If $A$ is Gorenstein and $M$ and $N$ are modules with $M \in
\mathcal{C}^{rc}_A$, then the implication  
$$\Ext_A^i(N,M)=0 \text{ for } i \gg 0 \hspace{2mm} \Rightarrow
\hspace{2mm} \Ext_A^i(M,N)=0 \text{ for } i \gg 0$$
holds. In particular, symmetry in the vanishing of cohomology holds
for modules having reducible complexity.
\end{proposition}

\section{Complexity testing}

In this section we introduce a method for computing an upper bound for
the complexity of a given module in $\mathcal{C}^{frc}_A$. We start
with a key result which shows that the modules in $\mathcal{C}^{frc}_A$
having infinite projective dimension also have higher self
extensions. The result is just a genralization of \cite[Corollary
3.2]{Bergh1} from $\mathcal{C}^{rc}_A$ to $\mathcal{C}^{frc}_A$, but
we include the proof for the convenience of the reader. 

\begin{proposition}\label{projdim}
If $M$ belongs to $\mathcal{C}^{frc}_A$, then
$$\Pd M = \sup \{ i | \Ext_A^i(M,M) \neq 0 \}.$$
\end{proposition}

\begin{proof}
If the projective dimension of $M$ is finite, then the first part of the
proof of \cite[Theorem 3.1]{Bergh1} shows that $\Ext_A^{\Pd M}(M,M)
\neq 0$. Suppose therefore that the projective dimension is
infinite. By definition, there exists a positive degree homogeneous
element $\eta \in \Ext_A^*(M,M)$ such that $\cx K_{\eta} = \cx
M-1$. Suppose now that $\Ext_A^i(M,M)$ vanishes for $i \gg 0$. Then 
$\eta$ is nilpotent, that is, there is a number $t$ such that $\eta^t
=0$. The exact sequence
$$0 \to M \to K_{\eta^t} \to \Omega_A^{t|\eta|-1}(M) \to 0$$
corresponding to $\eta^t$ then splits, and therefore $\cx K_{\eta^t}
= \cx M$ since the end terms are of the same
complexity. However, it follows from \cite[Lemma 2.3]{Bergh1} that for
any numbers $m$ and $n$ the $A$-modules $K_{\eta^m}$ and $K_{\eta^n}$
are related through an exact sequence
$$0 \to \Omega_A^{n|\eta|}(K_{\eta^m}) \to K_{\eta^{m+n}} \oplus F \to
K_{\eta^n} \to 0,$$
in which $F$ is some free module. Using this and the fact that in a
short exact sequence the complexity of the middle term is at most the
maximum of the complexitites of the end terms, an induction argument
gives the inequality $\cx K_{\eta^i} \le \cx K_{\eta}$ for every
$i \ge 1$. Combining all our obtained (in)equalities on complexity, we get
$$\cx M = \cx K_{\eta^t} \le  \cx K_{\eta} = \cx M-1,$$
a contradiciton. Therefore $\Ext_A^i(M,M)$ cannot vanish for all $i \gg 0$
when the projective dimension of $M$ is infinite.
\end{proof}

We are now ready to prove the main result. For a given natural number
$t$, denote by $\mathcal{C}^{frc}_A(t)$ the full subcategory 
$$\mathcal{C}^{frc}_A(t) \stackrel{\text{def}}{=} \{ X \in
\mathcal{C}^{frc}_A | \cx X =t \}$$
of $\mathcal{C}^{frc}_A$ consisiting of the modules of complexity
$t$. The main result shows that this subcategory serves as a
``complexity test category'', in the sense that if a module has no
higher extensions with the modules in $\mathcal{C}^{frc}_A(t)$, then
its complexity is strictly less than $t$.

\begin{theorem}\label{test}
Let $M$ be a module belonging to $\mathcal{C}^{frc}_A$ and $t$ a natural
number. If $\Ext_A^i(M,N)=0$ for every $N \in \mathcal{C}^{frc}_A(t)$
and $i \gg 0$, then $\cx M < t$.
\end{theorem}

\begin{proof}
We show by induction that if $\cx M \ge t$, then there is a module $N \in
\mathcal{C}^{frc}_A(t)$ with the property that $\Ext_A^i(M,N)$ does
not vanish for all $i \gg 0$. If the complexity of $M$ is $t$, then by
Proposition \ref{projdim} we may take $N$ to be $M$ itself, so suppose
that $\cx M > t$. Choose a cohomological homogeneous element $\eta \in
\Ext_A^*(M,M)$ of positive degree reducing the complexity. In the
corresponding exact sequence
$$0 \to M \to K_{\eta} \to \Omega_A^{|\eta|-1}(M) \to 0,$$
the module $K_{\eta}$ also belongs to $\mathcal{C}^{frc}_A$ and has
complexity one less than that of $M$, hence by induction there is a
module $N \in \mathcal{C}^{frc}_A(t)$ such that $\Ext_A^i(K_{\eta},N)$
does not vanish for all $i \gg 0$. From the long exact sequence
$$\cdots \to \Ext_A^{i+|\eta|-1}(M,N) \to
\Ext_A^i(K_{\eta},N) \to \Ext_A^i(M,N) \to \Ext_A^{i+|\eta|}(M,N) \to
\cdots$$ 
resulting from $\eta$, we see that $\Ext_A^i(M,N)$ cannot vanish for
all $i \gg 0$.   
\end{proof}

In particular, we can use the category $\mathcal{C}^{frc}_A(1)$ to
decide whether a given module in $\mathcal{C}^{frc}_A$ has finite
projective dimension. We record this fact in the following corollary.

\begin{corollary}\label{testcor}
A module $M \in \mathcal{C}^{frc}_A$ has finite projective
dimension if and only if $\Ext_A^i(M,N)=0$ for every $N \in
\mathcal{C}^{frc}_A(1)$ and $i \gg 0$.
\end{corollary}

\begin{remark}
Let $\mathcal{C}^{ci}_A$ denote the category of all $A$-modules of
finite complete intersection dimension, and for each natural number
$t$ define the two categories
\begin{eqnarray*}
\mathcal{C}^{ci}_A(t) & \stackrel{\text{def}}{=} & \{ X \in
\mathcal{C}^{ci}_A | \cx X =t \} \\
\mathcal{C}^{rc}_A(t) & \stackrel{\text{def}}{=} & \{ X \in
\mathcal{C}^{rc}_A | \cx X =t \}.
\end{eqnarray*}
Then Theorem \ref{test} and Corollary \ref{testcor} remain true if we
replace $\mathcal{C}^{frc}_A$ and $\mathcal{C}^{frc}_A(t)$ by
$\mathcal{C}^{rc}_A$ (respectively, $\mathcal{C}^{ci}_A$) and
$\mathcal{C}^{rc}_A(t)$ (respectively, $\mathcal{C}^{ci}_A(t)$). That
is, when the module we are considering has reducible complexity
(respectively, finite complete intersection dimension), then we need
only use modules of reducible complexity (respectively, finite
complete intersection dimension) as test modules. 
\end{remark}

When the ring is Gorenstein, then it follows from Proposition
\ref{symmetry} that symmetry holds for the vanishing of cohomology
between modules of reducible complexity. We therefore have the
following symmetric version of Theorem \ref{test}.   

\begin{corollary}\label{gorensteintest}
Suppose $A$ is Gorenstein, let $M$ be a module belonging to
$\mathcal{C}^{rc}_A$, and let $t$ be a natural number. If
$\Ext_A^i(N,M)=0$ for every $N \in \mathcal{C}^{rc}_A(t)$ and $i \gg
0$, then $\cx M < t$. 
\end{corollary}

We now turn to the setting in which \emph{every} $A$-module has
reducible complexity. For the remainder of this section, \emph{we assume $A$
is a complete intersection}, i.e.\  the $\m$-adic completion
$\widehat{A}$ of $A$ is the residue ring of a regular local ring
modulo a regular sequence. For such rings, Avramov and Buchweitz
introduced in \cite{Avramov1} and \cite{AvramovBuchweitz} a theory of
cohomological support varieties, and they showed that this theory is
similar to that of the cohomological support varieties for group
algebras. As we will implicitly use this powerful theory in the
results to come, we recall now the definitions (details can be found
in \cite[Section 1]{Avramov1} and \cite[Section 2]{AvramovBuchweitz}).  

Denote by $c$ the \emph{codimension} of
$A$, that is, the integer $\dim_k ( \m / \m^2 )- \dim A$, and by
$\chi$ the sequence $\chi_1, \dots, \chi_c$ consisting of 
the $c$ commuting Eisenbud operators of cohomological degree two. For
every $\widehat{A}$-module $X$ there is a homomorphism 
$$\widehat{A}[ \chi ] \xrightarrow{\phi_X} \Ext_{\widehat{A}}^*(X,X)$$
of graded rings, and via this homomorphism $\Ext_{\widehat{A}}^*(X,Y)$
is finitely generated over $\widehat{A}[ \chi ]$ for
every $\widehat{A}$-module $Y$. Denote by $H$ the polynomial ring $k [
\chi ]$, and by $E(X,Y)$ the graded space $\Ext_{\widehat{A}}^*(X,Y)
\otimes_{\widehat{A}} k$. The above homomorphism $\phi_X$, together with the
canonical isomorphism $H \simeq \widehat{A}[ \chi ]
\otimes_{\widehat{A}} k$, induce a homomorphism $H \to E(X,X)$ of
graded rings, under which $E(X,Y)$ is a finitely generated
$H$-module. Now let $M$ be an $A$-module, and denote by $\widehat{M}$
its $\m$-adic completion $\widehat{A} \otimes_A M$. The \emph{support
variety} $\V (M)$ of $M$ is the algebraic set
$$\V (M) \stackrel{\text{def}}{=} \{ \alpha \in \tilde{k}^c | f (
\alpha ) =0 \text{ for all } f \in \Ann_H E( \widehat{M}, \widehat{M}
) \},$$
where $\tilde{k}$ is the algebraic closure of $k$. Finally, for an
ideal $\az \subseteq H$ we define the variety $\V_H ( \az )
\subseteq \tilde{k}^c$ to be the zero set of $\az$. 

As mentioned above, this theory shares many properties with the theory
of cohomological support varieties for modules over group algebras of
finite groups. For instance, the dimension of the variety of a module
equals the complexity of the module, in particular the variety is
trivial if and only if the module has finite projective dimension. The
following complexity test result relies on \cite[Corollary
2.3]{Bergh2}, which says that every homogeneous algebraic subset of
$\tilde{k}^c$ is realizable as the support variety of some $A$-module. 

\begin{proposition}\label{vartest}
Let $M$ be an $A$-module, let $\eta_1, \dots, \eta_t \in H$ be
homogeneous elements of positive degrees, and choose an $A$-module
$T_{\eta_1, \dots, \eta_t}$ with the property that $\V ( T_{\eta_1,
  \dots, \eta_t}) = \V_H ( \eta_1, \dots, \eta_t )$. If
$\Ext_A^i(M,T_{\eta_1, \dots, \eta_t})=0$ for $i \gg 0$, then $\cx M
\le t$.   
\end{proposition}

\begin{proof}
Denote the ideal $\Ann_H E( \widehat{M}, \widehat{M} ) \subseteq H$ by 
$\az$. If $\Ext_A^i(M,T_{\eta_1, \dots, \eta_t})=0$ for $i \gg 0$,
then from \cite[Theorem 5.6]{AvramovBuchweitz} we obtain
$$\{ 0 \} = \V (M) \cap \V (T_{\eta_1, \dots, \eta_t}) = \V_H  ( \az )
\cap \V_H ( \eta_1, \dots, \eta_t ) = \V_H \left ( \az + ( \eta_1,
\dots, \eta_t ) \right ),$$
hence the ring $H / \left ( \az + ( \eta_1, \dots, \eta_t ) \right )$
is zero dimensional. But then the dimension of the ring $H / \az$ is
at most $t$, i.e.\ $\cx M \le t$. 
\end{proof}

We illustrate this last result with an example.

\begin{example}
Let $k$ be a field and $Q$ the formal power series ring $k \llbracket x_1,
\dots, x_c \rrbracket$ in $c$ variables. For each $1 \le i \le c$, let
$n_i \ge 2$ be an integer, let $\az \subseteq Q$ be the ideal
generated by the regular sequence $x_1^{n_1}, \dots, x_c^{n_c}$, and
denote by $A$ the complete intersection $Q / \az$. For each $1 \le i
\le c$ we shall construct an $A$-module whose support variety equals
$\V_H ( \chi_i )$, by adopting the techniques used in \cite[Section
7]{Snashall} to give an interpretation of the Eisenbud operators.

Consider the exact sequence
$$0 \to \m_Q \to Q \to k \to 0$$
of $Q$-modules. Applying $A \otimes_Q -$ to this sequence gives the
four term exact sequence
\begin{equation*}\label{ES}
0 \to \Tor^Q_1(A,k) \to \m_Q / \az \m_Q \to A \to k \to 0
\tag{\textdagger} 
\end{equation*}
of $A$-modules. Consider the first term in this sequence. By tensoring
the exact sequence
$$0 \to \az \to Q \to A \to 0$$
over $Q$ with $k$, we obtain the exact sequence
$$0 \to \Tor^Q_1(A,k) \to \az \otimes_Q k \xrightarrow{g} Q \otimes_Q
k \to A \otimes_Q k \to 0,$$
in which the map $g$ must be the zero map since $\az k =0$. This gives
isomorphisms 
$$\Tor^Q_1(A,k) \simeq \az \otimes_Q k \simeq \az \otimes_Q (A
\otimes_A k) \simeq \az / \az^2 \otimes_A k$$
of $A$-modules. Since $\az$ is generated by a regular sequence of
length $c$, the $A$-module $\az / \az^2$ is free of rank $c$, and
therefore $\Tor^Q_1(A,k)$ is isomorphic to $k^c$. We may now rewrite
the four term exact sequence (\ref{ES}) as
$$0 \to k^c \xrightarrow{f} \m_Q / \az \m_Q \to A \to k \to 0,$$
and it is not hard to show that the map $f$ is defined by
$$( \alpha_1, \dots, \alpha_c ) \mapsto \sum \alpha_i x_i^{n_i} + \az
\m_Q.$$ 
The image of the Eisenbud operator $\chi_j$ under the homomorphism
$\widehat{A} [ \chi ] \xrightarrow{\phi_k} \Ext_{\widehat{A}}^*(k,k)$
is the bottom row in the pushout diagram
$$\xymatrix{
0 \ar[r] & k^c \ar[r]^>>>>>f \ar[d]^{\pi_j} & \m_Q / \az \m_Q \ar[r]
\ar[d] & A \ar[r] \ar@{=}[d] & k \ar[r] \ar@{=}[d] & 0 \\
0 \ar[r] & k \ar[r] & K_{\chi_j} \ar[r] & A \ar[r] & k \ar[r] & 0}$$
of $A$-modules, in which the map $\pi_j$ is projection onto the $j$th
summand. The pushout module $K_{\chi_j}$ can be described explicitly
as
$$K_{\chi_j} = \frac{k \oplus \m_Q / \az \m_Q}{\{ ( \alpha_j,- \sum
  \alpha_i x_i^{n_i} + \az \m_Q ) \mid ( \alpha_1, \dots, \alpha_c )
  \in k^c \}},$$
and by \cite[Theorem 2.2]{Bergh2} its support variety is given by $\V
( K_{\chi_j} ) = \V (k) \cap \V_H ( \chi_j )$. But the variety of $k$
is the whole space, hence the equality $\V ( K_{\chi_j} ) = \V_H (
\chi_j )$. Thus by Proposition \ref{vartest} the $A$-module
$K_{\chi_j}$ is a test module for finding modules with bounded
projective resolutions; if $M$ is an $A$-module such that
$\Ext_A^i(M,K_{\chi_j})=0$ for $i \gg 0$, then $\cx M \le 1$.
\end{example}

Before proving the final result, we need a lemma showing that every
maximal  
Cohen-Macaulay module over a complete intersection has reducible
complexity by a cohomological element of degree two. This improves
\cite[Lemma 2.1(i)]{Bergh3}, which states that such a cohomological
element exists after passing to some suitable faithfully flat extension
of the ring. 

\begin{lemma}\label{reducing}
If $M$ is a maximal Cohen-Macaulay $A$-module of infinite projective
dimension, then there exists an element $\eta \in \Ext_A^2(M,M)$
reducing its complexity.
\end{lemma}

\begin{proof}
Since the dimension of $\V (M)$ is nonzero, the radical $\sqrt{\Ann_H E(
\widehat{M}, \widehat{M} )}$ of $\Ann_H E(
\widehat{M}, \widehat{M} )$ is properly contained in the graded maximal
ideal of $H$. Therefore one of the Eisenbud operators, say $\chi_j$,
is not contained in $\sqrt{\Ann_H E( \widehat{M}, \widehat{M}
  )}$. We now follow the arguments given prior to \cite[Corollary
2.3]{Bergh2}. Viewing $\chi_j$ as an element of $\widehat{A}[ \chi ]$, we can
apply the homomorphism $\phi_{\widehat{M}}$ and obtain the element
$\phi_{\widehat{M}} ( \chi_j ) \otimes 1$ in $\Ext_{\widehat{A}}^2(
\widehat{M}, \widehat{M} ) \otimes_{\widehat{A}} k$. Now $\Ext_{\widehat{A}}^2(
\widehat{M}, \widehat{M} )$ is isomorphic to $\Ext_A^2(M,M) \otimes_A
\widehat{A}$, and there is an isomorphism
$$\Ext_A^2(M,M) \otimes_A k \xrightarrow{\sim} \Ext_{\widehat{A}}^2(
\widehat{M}, \widehat{M} ) \otimes_{\widehat{A}} k$$
mapping an element $\theta \otimes 1 \in \Ext_A^2(M,M) \otimes_A k$ to
$\widehat{\theta} \otimes 1$. Therefore there exists an element $\eta
\in \Ext_A^2(M,M)$ such that $\widehat{\eta} \otimes 1$ equals
$\phi_{\widehat{M}} ( \chi_j ) \otimes 1$ in $\Ext_{\widehat{A}}^2( 
\widehat{M}, \widehat{M} ) \otimes_{\widehat{A}} k$. If the exact
sequence
$$0 \to M \to K_{\eta} \to \Omega_A^1 (M) \to 0$$
corresponds to $\eta$, then its completion
$$0 \to \widehat{M} \to \widehat{K_{\eta}} \to \Omega_{\widehat{A}}^1 (
\widehat{M} ) \to 0$$
corresponds to $\widehat{\eta}$, and so from \cite[Theorem
2.2]{Bergh2} we see that 
$$\V ( K_{\eta} ) = \V (M) \cap \V_H ( \chi_j ).$$
Since $\chi_j$ was chosen so that it ``cuts down'' the variety of $M$,
we must have $\dim \V ( K_{\eta} ) = \dim \V (M)-1$, i.e.\ $\cx
K_{\eta} = \cx M-1$.
\end{proof}

We have now arrived at the final result,
which improves Theorem \ref{test} when the ring is a complete
intersection. Namely, for such rings it suffices to check the
vanishing of finitely many cohomology groups ``separated'' by an odd
number. The number of cohomology groups we need to check depends on
the complexity value we are testing. Recall that we have denoted the
codimension of the complete intersection $A$ by $c$.

\begin{theorem}\label{testci}
Let $M$ be an $A$-module and $t \in \{ 1, \dots, c \}$ an integer. If
for every $A$-module $N$ of complexity $t$ there is an odd number $q$
such that 
$$\Ext_A^n(M,N) = \Ext_A^{n+q}(M,N) = \cdots = \Ext_A^{n+(c-t)q}(M,N) =0$$
for some even number $n> \dim A - \depth M$, then $\cx M <t$.
\end{theorem}   

\begin{proof}
Since $\cx M = \cx \Omega_A^{\dim A - \depth M}(M)$, we may without
loss of generality assume that $M$ is maximal Cohen-Macaulay and that 
$n>0$. We prove by induction that if $\cx M \ge t$, then for any odd 
number $q$ and any even integer $n >0$, the groups 
$$\Ext_A^n(M,N), \hspace{1mm} \Ext_A^{n+q}(M,N), \dots,
\Ext_A^{n+(\cx M-t)q}(M,N)$$   
cannot all vanish for every module $N$ of complexity $t$. When the
complexity of $M$ is $t$, take $N$ to be $M$ itself. In this case it
follows from \cite[Theorem 4.2]{AvramovBuchweitz} that $\Ext_A^n(M,N)$
is nonzero, because $t \ge 1$. Now assume $\cx M >t$, and write $q$ as
$q=2s-1$ where $s \ge 1$ is an integer. By Lemma \ref{reducing} there
is an element $\eta \in \Ext_A^2(M,M)$ reducing the complexity of $M$,
and it follows from \cite[Proposition 2.4(i)]{Bergh1} that the element 
$\eta^s \in \Ext_A^{2s}(M,M)$ also reduces the complexity. The latter
element corresponds to an exact sequence
$$0 \to M \to K \to \Omega_A^q(M) \to 0,$$
in which the complexity of $K$ is one less than that of $M$. By
induction there exists a module $N$, of complexity $t$, such that the
groups 
$$\Ext_A^n(K,N), \hspace{1mm} \Ext_A^{n+q}(K,N), \dots, \Ext_A^{n+(\cx
  K-t)q}(K,N)$$ 
do not all vanish. Then from the exact sequence we see
that the groups 
$$\Ext_A^n(M,N), \hspace{1mm} \Ext_A^{n+q}(M,N), \dots, \Ext_A^{n+(\cx
  M-t)q}(M,N)$$ 
cannot possible all vanish. Since the complexity of any $A$-module is
at most $c$, the proof is complete.
\end{proof}

\end{document}